# Regular and oscillatory parts
# for basic functions of prime numbers
# I. Regular parts

R. M. Abrarov[1] and S. M. Abrarov[2]

## Abstract

We consider the regular parts for basic functions of prime numbers with Riemann approximation accuracy.

## 1. Introduction

Apparently the importance of the distribution of prime numbers was first recognized in connection between prime numbers and the zeta function $\zeta(s)$ by the celebrated Euler product

$$\zeta(s) \equiv \sum_{n \geq 1} \frac{1}{n^s} = \frac{1}{\prod_{p(primes)}\left(1-\frac{1}{p^s}\right)}, \qquad (1)$$

published in Euler's book "Introductio in Analysin Infinitorum" in 1748. This identity shows us not only the infinitude of the prime numbers

$$\pi(x) \equiv \sum_{p \leq x} 1 \to \infty \ , \text{ when } x \to \infty, \qquad (2)$$

but also gives us the comprehension about the distribution of primes.

Legendre in the second volume of his book "Number Theory", published in 1808 [1], proposed an approximate formula enumerating primes up to a given value *x*. Analyzing available prime tables, Legendre came to conclusion that the best approximation for the distribution of primes may be given by

$$\pi(x) \approx \frac{x}{\log x - A}, \qquad (3)$$

with $A = 1.08366$. Formula (3) seemed to be in a good agreement with prime table for the interval from 10,000 up to 1,000,000.

Gauss demonstrated rather different approach to this problem. As we can see from the letter to his apprentice, astronomer Hencke, elderly Gauss in 1849 wrote that in his early years he empirically discovered that "the density of primes around $x$ is about $1/\log x$" [2] and derived from this observation an excellent approximation to $\pi(x)$ through the logarithmic integral function $Li(x)$:

$$\pi(x) \approx Li(x) \equiv \int_2^x \frac{dt}{\log t}. \qquad (4)$$

Such an original probabilistic approach had further development in the works of Cramer and his successors [2, 3].

Chebyshev [tʃebiˈʃɔːv] made first analytical contribution to this problem in 1849 in his memoir "On the definition of number of primes not exceeding a given magnitude" [4]. Using Euler product (1) with real variable $s$, he proved that Legendre conjecture can be true only with constant A equal to 1. Moreover he showed that $Li(x)$- function gives even a better estimation for prime counting function $\pi(x)$.

In the next memoir "On the prime numbers" (1852) he introduced two new important counting functions (now known as Chebyshev's functions) $\theta(x) = \sum_{p \leq x} \log p$ and $\psi(x) = \sum_{p^n \leq x} \log p$, which have the following properties at $x \to \infty$: if any of three limits $\frac{\psi(x)}{x}$, $\frac{\theta(x)}{x}$, $\frac{\pi(x)\log x}{x}$ exist then all of them must exist and be equal to 1. However Chebyshev did not provide any proof that this limit exists.

The statement

$$\pi(x) \sim \frac{x}{\log x} \quad (x \to \infty) \qquad (5)$$

is the subject of the Prime Number Theorem, which was independently proved by Hadamard [5] and de la Vallée Poussin [6] in 1896. This relation does not contradict to Lejendre and Gauss estimation since expressions (3), (4) and (5) have the same limit.

43The approach of Riemann [1, 7] to this problem showed that $\pi(x)$-function can be approximated by

$$Ri(x) = \sum_{n=1}^{\infty} \frac{\mu(n)}{n} Li(x^{1/n}), \qquad (6)$$

where $\mu(n)$ is Möbius function defined as

$$\mu(n) = \begin{cases} 1 & \text{if } n = 1, \\ (-1)^m & \text{if } n \text{ is a product of } m \text{ distinct primes,} \\ 0 & \text{if the square of primes divides } n. \end{cases} \qquad (7)$$

Equation (6) is known as Riemann's approximation to $\pi(x)$. In comparison with *Li(x)*-function, the formula (6) provides much better approximation at least in the range where the exact prime counting value is feasible for computing.

All given above estimations for $\pi(x)$-function describe only smooth regular behavior of this function and do not provide any information about irregularities (oscillations).

Riemann in his innovative approach demonstrated [7] that irregularities in the distribution of primes are closely related to the location of zeros of $\zeta(s)$-function (1) with the complex variable *s*. In this work Riemann assumed a famous hypothesis that all zeroes of $\zeta(s)$-function except the trivial ones lie on the *critical line Re(s) = ½*. This hypothesis is currently considered as one of the most celebrated open problems in mathematics [8, 9].

Von Koch showed [10] that the Riemann hypothesis is equivalent to

$$\pi(x) = Li(x) + O(\sqrt{x} \log x), \qquad (8a)$$

or, in more general form, to

$$\pi(x) = Li(x) + O_\varepsilon(x^{1/2+\varepsilon}) \qquad (\forall \varepsilon > 0), \qquad (8b)$$

where $\varepsilon$ means that in the expression for $O_\varepsilon$ the argument may be product of $x^{1/2}$ and any fixed power of *log(x)*.

The brief historical overview and some facts, described above may be helpful in reading this article. However, for more complete information the interested readers are referred to the literature for state-of-the-art [3, 11-13].

In the article, the expressions for regular and oscillatory parts for basic functions of prime numbers are considered. The structure of the present work consists of two parts. In the first part, expressions for regular part of basic functions ($\theta(x), \psi(x), R(x)$) with Riemann approximation accuracy are shown. The second part, which will be published soon, describes the expressions for oscillatory parts of these functions.

## 2. Regular parts for $\pi(x)$, $R(x)$, $\theta(x)$, $\psi(x)$

Let us represent $\pi(x)$ through its regular $\pi_{reg}(x)$ and oscillatory parts $\pi_{osc}(x)$, where as a regular part of $\pi_{reg}(x)$ we imply Riemann approximation (6)

$$\pi(x) = \pi_{reg}(x) + \pi_{osc}(x) = \sum_{n=1}^{\infty} \frac{\mu(n)}{n} Li(x^{1/n}) + \pi_{osc}(x). \tag{9}$$

In order to derive an identity between $\pi(x)$ and $Li(x)$ functions, define

$$R(x) \equiv \sum_{p \leq x} \frac{\log p}{p} - \log x, \tag{10}$$

where $p$ runs over all primes up to $x$.

For convenience we also introduce complementary function $r(x) = \sum_{p \leq x} \frac{\log p}{p}$. Then using Stietjes integration method and integrating by parts, we can write the following

$$\pi(x) = \sum_{k=2}^{x}(r_k - r_{k-1})\frac{k}{\log k} = \int_{2-\varepsilon}^{x} \frac{y}{\log y} dr(y) = r(x)\frac{x}{\log x} - \int_{2}^{x} r(y) d\left(\frac{y}{\log y}\right). \tag{11}$$

Consequently, rewriting (11) in terms of $R(x)$, we obtain the desired identity [14]

$$\pi(x) = Li(x) + R(x)\frac{x}{\log x} - \int_{2}^{x} R(y) d\left(\frac{y}{\log y}\right) + 2, \tag{12}$$

where[*]

$$Li(x) \equiv \int_{2}^{x} \frac{1}{\log y} dy. \tag{13}$$

Taking derivative of this identity, we have

$$\pi(x)' = \frac{1}{\log x} + \frac{xR(x)'}{\log x} \tag{14}$$

___

[*] It should be noted that the definition of the logarithmic integral is the matter of our convenience and therefore it can be represented, for instance, in the form of the "Cauchy principal value" [13]:

$$Li(x) \equiv \lim_{\varepsilon \to 0+}\left(\int_{0}^{1-\varepsilon} + \int_{1+\varepsilon}^{x}\right)\frac{1}{\log t} dt.$$

This definition of logarithmic integral differs form that of (13) just by constant ≈ 1.04. In our recent publication, the same result (12) has been obtained by Dirac delta function approach [14], where both definitions of logarithmic integral can be implied.

with corresponding counterpart

$$R(x)' = \frac{\pi(x)' \log x}{x} - \frac{1}{x}. \tag{15}$$

Recall Ramanujan formula, which actually is the direct sequence of Riemann approximation for $\pi_{reg}(x)$:

$$\frac{d\pi_{reg}(x)}{dx} = \frac{1}{x \log x} \sum_{k=1}^{\infty} \frac{\mu(k)}{k} x^{1/k}. \tag{16}$$

Substituting Ramanujan formula into (15), we get for regular part of $R(x)$

$$\frac{dR_{reg}(x)}{dx} = \sum_{k=2}^{\infty} \frac{\mu(k)}{k} x^{1/k-2}. \tag{17}$$

Integrating this expression, we obtain the formula

$$R_{reg}(x) = \frac{1}{x} \sum_{k=2}^{\infty} \frac{\mu(k)}{1-k} x^{1/k} + const. \tag{18}$$

It follows from the Prime Number Theorem that $R(x)$-function has to have the limit, which can be represented as [15, 16]

$$\lim_{x \to \infty} R(x) = -\gamma - \sum_{k=2}^{\infty} \sum_{p} \frac{\log p}{p^k} = -1.33258227573322087.... \tag{19}$$

Hence we see that in (18) $const. = -\gamma - \sum_{k=2}^{\infty} \sum_{p} \frac{\log p}{p^k}$.

For arbitrary $x$, $R(x)$ can be represented as

$$R(x) = R_{reg}(x) + R_{osc}(x)$$

with oscillatory part $R_{osc}(x)$, where for the regular part $R_{reg}(x)$ we have

$$R_{reg}(x) = -\gamma - \sum_{k=2}^{\infty} \sum_{p} \frac{\log p}{p^k} + \frac{1}{x} \sum_{k=2}^{\infty} \frac{\mu(k)}{1-k} x^{1/k} \approx -\gamma - \sum_{k \geq 2} \sum_{p^k \leq x} \frac{\log p}{p^k}. \tag{20}$$

Now we can complete the discussion concerning $R(x)$ function and let us proceed our consideration to Chebyshev's function

$$\theta(x) = \sum_{p \leq x} \log p.$$

Applying again the same procedure as for (11), we have

$$\pi(x) = \sum_{k=2}^{x} (\theta_k - \theta_{k-1}) \frac{1}{\log k} = \int_{2-\varepsilon}^{x} \frac{d\theta(y)}{\log y} = \frac{\theta(x)}{\log x} - \int_{2}^{x} \theta(y) d\left(\frac{1}{\log y}\right) \tag{21}$$

Differentiating gives us

$$\pi(x)' = \frac{\theta(x)'}{\log x}, \tag{22}$$

with its counterpart

$$\theta(x)' = \pi(x)' \log x. \tag{23}$$

Again applying Ramanujan formula (16) we have

$$\frac{d\theta_{reg}(x)}{dx} = \frac{1}{x}\sum_{k=1}^{\infty}\frac{\mu(k)}{k}x^{1/k}. \tag{24}$$

Integrating (24) results to

$$\theta_{reg}(x) = \sum_{k=1}^{\infty}\mu(k)x^{1/k} + const. \tag{25}$$

The constant in (25) is determined through

$$\theta_{reg}(0) = const. = 0. \tag{26}$$

Thus $\theta(x)$ can be determined through following equation

$$\theta(x) = \sum_{k=1}^{\infty}\mu(k)x^{1/k} + \theta_{osc}(x), \tag{27}$$

where $\theta_{osc}(x)$ is the oscillatory part of $\theta(x)$.

Recall the formula relating two Chebyshev functions $\theta(x)$ and $\psi(x)$ [13]

$$\theta(x) = \sum_{k \geq 1}\mu(k)\psi\left(x^{1/k}\right). \tag{28}$$

Hence we see that the regular part of Chebyshev's function $\psi(x)$ can be expressed as

$$\psi_{reg}(x) = x. \tag{29}$$

As a result, $\psi(x)$ function can be represented as follows

$$\psi(x) = x + \psi_{osc}(x), \tag{30}$$

where $\psi_{osc}(x)$ is oscillatory part of this function.

Now we can see that Riemann's approximation (9) inevitably leads to (20), (27) and (29) for regular parts of $R(x)$, $\theta(x)$ and $\psi(x)$, respectively. Contrarily, assumption (29) leads to Riemann's approximation (9) for the regular part of $\pi(x)$-function.

## 3. Conclusion

The counting functions $\pi(x)$, $R(x)$, $\theta(x)$ and $\psi(x)$ can be expressed explicitly through corresponding regular and oscillatory parts. Using Riemann's approximation for $\pi(x)$-function, the regular parts of $R(x)$, $\theta(x)$ and $\psi(x)$ functions were obtained. Contrarily, an assumption (29) for the regular part of $\psi(x)$-function leads to Riemann's approximation (9).

## Acknowledgements

Authors express gratitude to Dr. J. E. Stirling and Dr. E. N. Tsoy for support and fruitful discussions.

[1] University of Toronto, Canada           rabrarov@physics.utoronto.ca

[2] York University, Toronto, Canada        abrarov@yorku.ca
Dongguk University, Seoul, South Korea      abrarov@dongguk.edu